\documentclass{article}
\usepackage{amssymb}
\usepackage{graphicx}
\usepackage{amsmath}

\newtheorem{theorem}{Theorem}

\newtheorem{corollary}[theorem]{Corollary}

\newtheorem{lemma}[theorem]{Lemma}

\begin{document}

\title{Group laws and free subgroups in topological groups}
\author{Mikl\'{o}s Ab\'{e}rt}
\date{June 25, 2003}
\maketitle

\begin{abstract}
We prove that a permutation group in which different finite sets have
different stabilizers cannot satisfy any group law. For locally compact
topological groups with this property we show that almost all finite subsets
of the group generate free subgroups.

We derive consequences of these theorems on Thompson's group $F$, weakly
branch groups, automorphism groups of regular trees and profinite groups
with alternating composition factors of unbounded degree.
\end{abstract}

\section{Introduction}

Let $G$ be a permutation group acting on a set $X$. We say that $G$ \emph{%
separates} $X$ if for any finite subset $Y\subseteq X$ the pointwise
stabilizer $G_{Y}$ does not stabilize any point outside $Y$. A group $G$
satisfies a group law, if there exists a word $w\in F_{k}$ such that $%
w(g_{1},g_{2},\ldots ,g_{k})=1$ for all $g_{i}\in G$. Our first result is
the following.

\begin{theorem}
\label{nolaw}If $G$ separates $X$ then $G$ does not satisfy any group law.
\end{theorem}

In contrast, for the natural transitive permutation action of the so-called
lamplighter group $L$ (the wreath product $C_{2}\wr \mathbb{Z}$) the
stabilizer of $n$ points stabilizes at most $n$ other points. On the other
hand, $L$ is metabelian, so it satisfies a commutator identity.

Theorem \ref{nolaw} allows us to give a short proof of the main result in 
\cite{brin}. Thompson's group $F$ is the group of piecewise linear
orientation-preserving homeomorphisms of the unit interval with dyadic
breaking points and slopes which are powers of $2$.

\begin{corollary}
\label{thompson}Thompson's group $F$ does not satisfy any group law.
\end{corollary}

Note that by \cite{brin} $F$ does not contain any free subgroups.

We can generalise Theorem \ref{nolaw} for locally compact topological groups
which act on a set separating it and show that generic subgroups of them are
free. Altough we derive the compact case as a simple consequence of the
general result on locally compact groups, we state it differently as the
necessary restrictions are more natural.

Let $G$ be a compact topological group. Then there is a unique $G$-invariant
probability measure $\mu $ on $G$. This extends to a probability measure on $%
G^{n}$. A subset of $G^{n}$ is almost sure (a.s.) if its complement has zero
probability. An action of $G$ on a set $X$ is \emph{topological}, if for all 
$x\in X$ the point stabilizer $G_{x}$ is closed. If $G$ is also transitive,
then $X$ can be identified with the coset space $G/H$ for any point
stabilizer $H$ and as such, it carries a natural $G$-invariant measure $%
\lambda $.

A group $\Gamma \subseteq Sym(\Omega )$ acts freely on $\Omega $ if for
every $x\in \Omega $ the point-stabilizer $\Gamma_{x}$ is trivial. A group $%
H\subseteq G$ acts \emph{almost freely} on $X$ if there is $Y\subseteq X$
such that $\lambda (X\backslash Y)=0$, $Y$ is $\Gamma$-invariant and $\Gamma$
acts freely on $Y$.

\begin{theorem}
\label{compactfree}Let $G$ be a compact topological group acting
topologically on $X$ and separating $X$. Then $n$ random elements in $G$
generate a free subgroup of rank $n$ a.s. If $G$ is transitive then this
group acts almost freely on $X$ a.s.
\end{theorem}

So somewhat surprisingly, in the transitive case a generic finitely
generated subgroup of $G$ does exactly the opposite of separating $X$, as
almost all points have a stabilizer which stabilizes almost all points in $X$%
.

Theorem \ref{compactfree} applies to profinite weakly branch groups. A
weakly branch group can be defined as a group acting spherically
transitively on a rooted tree such that the rigid stabilizer of every vertex
is nontrivial. Recall that the rigid stabilizer of a vertex $v$ is the set
of group elements which move only descendants of $v$. The class of weakly
branch groups contains many groups with interesting properties. Examples are
the first Grigorchuk group (\cite{grigor}), Wilson's group of non-uniform
exponential growth (\cite{wilnonuni}), Segal's group of slow subgroup growth
(\cite{segal}) and the group constructed by Grigorchuk and Zuk which is not
in the class SG (\cite{grizuk}) but as Bartholdi and Vir\'{a}g have recently
shown is amenable (\cite{barvir}).

A profinite weakly branch group is a weakly branch group which is closed in
the profinite topology of the automorphism group of the corresponding tree,
e.g., the closure of any weakly branch group.

\begin{corollary}
\label{branch}If $G$ is a profinite weakly branch group then $n$ random
elements generate a free subgroup of rank $n$ in $G$ a.s. Moreover, if the
underlying tree has bounded degree and $G$ is finitely generated then $G$
contains a finitely generated dense free subgroup.
\end{corollary}

The above result generalises a result of Wilson \cite{wilson} who proved
that for branch profinite groups the subset of $n$-tuples in $G^{n}$
generating non-free subgroups is meagre (first category) in $G^{n}$. His
result implies that finitely generated pro-$p$ branch groups have finitely
generated dense free subgroups, but for the general profinite setting our
stronger probabilistic result is needed.

Theorem \ref{compactfree} also applies to infinite iterated wreath products
of nontrivial finite permutation groups. This reproves a result of
Bhattacharjee \cite{bhatta} saying that in these groups $n$ elements
generate a free subgroup of rank $n$ a.s.

We deduce Theorem \ref{compactfree} as a special case of a more general
theorem on locally compact groups. If the group $G$ is locally compact but
not compact, we do not have a $G$-invariant probability measure anymore.
Also, we have to strengthen the notion of topological action, since a
countable group can easily act on a set, separating it, without containing
any free subgroups. Thompson's group or the finitary symmetric group acting
on $\mathbb{N}$ both have this property.

A locally compact group $G$ carries a right $G$-invariant Haar measure $\mu $%
. An action of $G$ on a set $X$ is topological, if for all $x\in X$ the
point stabilizer $G_{x}$ is closed and $\mu (G_{x})=0$. The action is \emph{%
strongly topological}, if for all finite subset $Y\subseteq X$ the pointwise
stabilizer $G_{Y}$ acts topologically on $X\backslash Y$. A property of $n$%
-tuples of $G$ is almost sure (a.s.) if its complement is $\mu $-negligible
in $G^{n}$. Note that these notions are independent of the choice of $\mu $.
If $G$ is transitive, then $X$ again can be identified with the coset space $%
G/H$ for some closed subgroup $H$ and as such, it inherits a natural measure 
$\lambda $ from $G$. This measure is not necessarily $G$-invariant, but
always $G$-equivariant.

\begin{theorem}
\label{topfree}Let $G$ be a locally compact topological group acting
strongly topologically on $X$ and separating $X$. Then $n$ elements in $G$
generate a free subgroup of rank $n$ a.s. If $G$ is transitive then this
group acts almost freely on $X$ a.s.
\end{theorem}

A natural class of locally compact groups to which Theorem \ref{topfree}
applies are automorphism groups of unrooted regular trees.

\begin{corollary}
\label{regtree}For any regular tree $T$, $n$ elements in $Aut(T)$ generate a
free subgroup of rank $n$ a.s.
\end{corollary}

This corollary is used in \cite{abgla} to show that $Aut(T)$ contains
finitely generated dense free subgroups for every unrooted regular tree $T$.

Using an asymptotic version of Theorem \ref{compactfree} we can also answer
a question of Pyber. Dixon, Pyber, Seress and Shalev in \cite{pysha} gave an
elegant new proof to a conjecture of Magnus. The original proof was given by
Weigel \cite{weigel} using different methods. Magnus conjectured that the
free group $F_{2}$ is residually $S$ for every infinite set $S$ of pairwise
non-isomorphic nonabelian finite simple groups. Recall that a group $G$ is
residually $S$ if the intersection of all normal subgroups with quotients in 
$S$ is trivial in $G$. One of the corollaries in \cite{pysha} is that if a
profinite group $G$ has infinitely many nonabelian simple quotients then $n$
random elements in $G$ generate a free subgroup a.s. for all $n$.

It is asked in \cite[Problem 7]{pysha} whether one can generalise this to
profinite groups with infinitely many nonabelian simple composition factors.
Our result makes a first step into this direction.

\begin{theorem}
\label{alter}Let $G$ be a profinite group which has alternating composition
factors of unbounded degree. Then $n$ random elements generate a free
subgroup in $G$ a.s. for all $n$.
\end{theorem}

\section{Proofs}

In this section we prove the theorems and corollaries from the introduction.

\bigskip

\textbf{Proof of Theorem \ref{nolaw}.} We can assume that $G$ acts
transitively on $X$. If the action is intransitive, let $X^{\prime }$ be a
nontrivial orbit of $G$. Then $X^{\prime }$ is infinite, otherwise for every 
$x^{\prime }\in X^{\prime }$ and $Y=X^{\prime }\backslash \{x^{\prime }\}$
the stabilizer $G_{Y}$ would stabilize $x^{\prime }$, contradicting the
separability. By the same argument for every finite $Y\subseteq X$ the
orbits of the action of $G_{Y}$ on $X\backslash Y$ are all infinite.

Let $F_{k}$ be the free group on $k$ variables $f_{1},f_{2},\ldots ,f_{k}$.
Let $w\in F_{k}$ be a reduced word of length $n$, i.e., 
\begin{equation*}
w=v_{1}v_{2}\cdots v_{n}
\end{equation*}
where 
\begin{equation*}
v_{i}\in \left\{ f_{1},f_{1}^{-1},\ldots ,f_{k},f_{k}^{-1}\right\}
\end{equation*}
with $v_{i}v_{i+1}\neq 1$ for $1\leq i\leq n-1$.

Let $w_{j}$ be the $j$-th beginning segment of $w$, 
\begin{equation*}
w_{j}=v_{1}v_{2}\cdots v_{j}\text{ \ }(0\leq j\leq n)\text{.}
\end{equation*}

For a $k$-tuple $(g_{1},g_{2},\ldots ,g_{k})\in G^{k}$ we denote $%
w_{j}(g_{1},g_{2},\ldots ,g_{k})\in G$ by the value of the word $w_{j}$ in $%
G $ via the substitution $f_{i}=g_{i}$.

We claim that for all $1\neq w\in F_{k}$ and for all $x_{0}\in X$ there
exist a $k$-tuple $(h_{1},h_{2},\ldots ,h_{k})\in G^{k}$ such that the
points 
\begin{equation*}
x_{j}=x_{0}^{w_{j}(h_{1},h_{2},\ldots ,h_{k})}\text{ \ }(0\leq j\leq n)
\end{equation*}
are disjoint. This clearly shows that $G$ cannot satisfy the group law $w$,
as this particular substitution moves the point $x_{0}$ to $x_{n}\neq x_{0}$.

We prove our claim by induction on $n$. For $n=1$ the word consists of one
variable and the claim is trivial.

Fix $x_{0}\in X$. Using the inductional hypothesis, we have $%
(g_{1},g_{2},\ldots ,g_{k})\in G^{k}$ such that $x_{0},x_{1},\ldots ,x_{n-1}$
are all disjoint. If $x_{n}\notin \left\{ x_{0},x_{1},\ldots
,x_{n-1}\right\} $ then we found the right $(g_{1},g_{2},\ldots ,g_{k})$. So
assume that $x_{n}=x_{j}$ for some $j<n$.

Now let 
\begin{equation*}
I=\left\{ i<n\mid v_{i}=v_{n}\text{ or }v_{i+1}=v_{n}^{-1}\right\} \text{.}
\end{equation*}

Let $m$ be the index of $v_{n}$, that is $v_{n}\in \left\{
f_{m},f_{m}^{-1}\right\} $. Also let $g=v_{n}(g_{1},g_{2},\ldots ,g_{k})$
(thus $g$ is either $g_{m}$ or $g_{m}^{-1}$). We claim that $j\notin I$. If $%
v_{j}=v_{n}$, then $x_{j-1}^{g}=x_{j}=x_{n-1}^{g}$ which is a contradiction
since $j<n$ and all the $x_{i}$ are disjoint for $i<n$. If $%
v_{j+1}=v_{n}^{-1}$ then $x_{j+1}^{g}=x_{j}=x_{n-1}^{g}$ which implies $%
j=n-2 $ and thus $v_{n-1}=v_{n}^{-1}$ which contradicts the assumption that
the word $w$ is reduced. So our claim holds and $j\notin I$.

Now let $Y=\left\{ x_{i}\mid i\in I\right\} $ and let $c\in G_{Y}$ be an
arbitrary element of the stabilizer of $Y$.

Modify the $k$-tuple $(g_{1},g_{2},\ldots ,g_{k})$ to $(h_{1},h_{2},\ldots
,h_{k})$ as follows: $h_{i}=g_{j}$ ($i\neq m$), $h_{m}=g_{m}c$ if $%
v_{n}=f_{m}$ and $h_{m}=cg_{m}$ if $v_{n}=f_{m}^{-1}$. Using the
construction of $I$ it is straightforward to check that $%
x_{i}^{h_{i+1}}=x_{i+1}$ for $i<n-1$ and $x_{n-1}^{h_{n}}=x_{j}^{c}$.

Since $x_{j}\notin Y$, the orbit of $G_{Y}$ containing $x_{j}$ is infinite
so we can choose $c$ such that $x_{j}^{c}\notin \left\{ x_{0},x_{1},\ldots
,x_{n-1}\right\} $. So the $k$-tuple $(h_{1},h_{2},\ldots ,h_{k})$ shows
that our claim holds and so $w(g_{1},g_{2},\ldots ,g_{k}) \neq 1$. $\Box $

\bigskip

\textbf{Proof of Corollary \ref{thompson}. }We show that Thompson's group $F$
separates the unit interval by its natural action. Let $Y\subseteq \lbrack
0,1]$ be a finite set and let $x\in \lbrack 0,1]\backslash Y$. Let $%
d_{1},d_{2}\in (0,1)$ be dyadic numbers such that $x\in \lbrack d_{1},d_{2}]$
but $[d_{1},d_{2}]\cap Y=\emptyset $. Then the pointwise stabilizer $%
F_{d_{1},d_{2}}$ of $[0,d_{1}]\cup \lbrack d_{2},1]$ in $F$ is permutation
isomorphic to $F$ by its action on $[d_{1},d_{2}]$. In particular, there is
an element of $F_{d_{1},d_{2}}$ which moves the point $x$. $\Box $

\bigskip

\textbf{Proof of Theorem \ref{compactfree}. } We need to check that the
conditions for Theorem \ref{topfree} hold, that is, the action of $G$ on $X$
is strongly topological.

By the same argument as in the proof of Theorem \ref{nolaw} we see that
every $G$-orbit has to be infinite. So any point-stabilizer $G_{x}$ must
have infinite index in $G$; since $G_{x}$ is closed and $G$ is compact, $%
G_{x}$ has zero Haar measure, so the action of $G$ is topological. Also, for
any finite set $Y\subseteq X$ the pointwise stabilizer $G_{Y}$ trivially
separates $X\backslash Y$ and in $X\backslash Y$ the point-stabilizers are
closed (being the intersection of closed subgroups). The stabilizers $G_{Y}$
are compact, so by repeating the above argument to the action of $G_{Y}$ on $%
X\backslash Y$ we see that the action of $G$ on $X$ is strongly topological. 
$\Box $

\bigskip

\textbf{Proof of Corollary \ref{branch}. }Let $G$ act on the rooted tree $T$
spherically transitively such that the rigid vertex stabilizers are
nontrivial. Let $X$ be the boundary of the tree, i.e., the set of infinite
rays of $T$. Since $G$ is closed in the profinite topology of $Aut(T)$, for
every $x\in X$ the stabilizers $G_{x}$ are closed, so $G$ acts on $X$
topologically.

We claim that $G$ separates $X$. Let $y_{1},\ldots y_{n},x\in X$ be pairwise
disjoint rays. Choose a level $k$ such that the vertices in the rays $%
y_{1},\ldots y_{n},x$ at the $k$-th level are all disjoint. Let $x_{0}$ be
the vertex of $x$ at level $k$. Let $S$ be the stabilizer of $x_{0}$ in $G$
and let $R$ be the rigid vertex stabilizer of $x_{0}$ in $G$. Then $S$ acts
spherically transitively on the infinite subtree rooted at $x_{0}$, since if 
$a$ and $b$ are both descendants of $x_{0}$ at the same level, then there is 
$g\in G$ such that $a^{g}=b$ and clearly $g$ must stabilize $x_{0}$. Now $R$
is normal in $S$ and since it is nontrivial, it cannot stabilize any
infinite ray going through $x_{0}$. In particular, there exists $r\in R$
such that $x^{r}\neq x$. On the other hand, $r$ stabilizes every ray not
going through $x_{0}$. This proves our claim.

Now using Theorem \ref{compactfree} we see that $n$ random elements generate
a free subgroup of $G$ with probability $1$.

If the degree of $T$ is bounded, then $Aut(T)$ has only finitely many
different simple upper composition factors. Thus the degree of upper
alternating sections of $G$ is bounded. Recall that an upper section is a
homomorphic image of an open subgroup. Using \cite{bor} this implies that $G$
has polynomial maximal subgroup growth (PMSG), i.e., the number of maximal
subgroups of index $n$ in $G$ is bounded by a polynomial of $n$. By \cite
{mansha} for the profinite group $G$ PMSG is equivalent to that $G$ is
positively finitely generated (PFG), i.e., that there is a number $k$ such
that $k$ random elements generate a dense subgroup of $G$ with positive
probability.

So for every $n\geq k$ we obtain that $n$ random elements of $G$ generate a
dense free subgroup of $G$ with positive probability. This completes the
proof. $\Box $

\bigskip

Now we will prove Theorem \ref{topfree}. We need a straightforward lemma,
the proof of which we include for completeness.

If $G$ is a locally compact group and $H$ is a closed subgroup endowed with
a right Haar measure $\lambda $, then $\lambda $ naturally extends to every
right $H$-coset $C$ by $\lambda (X)=\lambda (Xg^{-1})$ where $g\in C$. The
measure $\lambda (X)$ is independent of the choice of $g$ as $\lambda $ is
right invariant in $H$.

Let $K\subset H$ be a closed subgroup of measure zero and let $\mu $ be
right Haar measure on $K$.

\bigskip

\begin{lemma}
\label{izelemma}Let $X$ be a $\lambda $-measurable subset of $C$ such that
for almost all $K$-cosets $D\subseteq C$ we have $\mu (D\cap X)=0$. Then $%
\lambda (X)=0$.
\end{lemma}

\textbf{Proof. }Fix $c_{0}\in C$. Then $C=Hc_{0}$. Define the subset $%
khc_{0} $%
\begin{equation*}
Y=\left\{ (k,h)\in K\times H\mid khc_{0}\in X\right\} \text{.}
\end{equation*}

Since $Y$ is the preimage of $X$ by a continuous map from $K\times H$ to $C$%
, it is measurable. Now for all $h_{0}\in H$ the fiber 
\begin{equation*}
R(h_{0})=\left\{ k\in K\mid (k,h_{0})\in Y\right\} =K\cap
Xc_{0}^{-1}h_{0}^{-1}
\end{equation*}
and so by the assumption for almost all $h_{0}$ we have 
\begin{equation*}
\mu (R(h_{0}))=\mu (Kh_{0}c_{0}\cap X)=0\text{.}
\end{equation*}

So by Fubini's theorem $Y$ has measure zero in $K\times H$. This implies
that for almost all $k_{0}\in K$ the fiber 
\begin{equation*}
L(k_{0})=\left\{ h\in H\mid (k_{0},h)\in Y\right\} =k_{0}^{-1}Xc_{0}^{-1}
\end{equation*}
has measure zero in $H$. Pick such a $k_{0}$. Using that $\lambda $ is left $%
H$-equivariant we get $\lambda (Xc_{0}^{-1})=0$ in $H$, that is, $\lambda
(X)=0$. $\Box $

\bigskip

\textbf{Proof of Theorem \ref{topfree}. }We will use the notations from the
proof of Theorem \ref{nolaw}.

We claim that we can assume that $G$ is transitive on $X$. Let $%
X_{0}\subseteq X$ be a $G$-orbit. The restriction $\overline{G}$ of the
action of $G$ on $X_{0}$ is a group homomorphism and the kernel $K$ is
closed of Haar measure zero since it is the intersection of the closed
stabilizers of points in $X_{0}$. Now $\overline{G}$ clearly separates $%
X_{0} $ and it acts strongly topologically on it. Since $K$ has zero
measure, the preimage of zero measure subsets of $\overline{G}^{n}$ have
zero measure in $G^{n}$ and our claim holds.

Let $Y_{1},Y_{2},\ldots ,Y_{k}\subseteq X$ and $Z_{1},Z_{2},\ldots
,Z_{k}\subseteq X$ be finite ordered subsets such that for all $1\leq i\leq
k $ we have $\left| Y_{i}\right| =\left| Z_{i}\right| $ and there exists $%
h_{i} $ such that $Y_{i}^{h_{i}}=Z_{i}$ as ordered subset. Let $Y=\cup Y_{i}$%
, $Z=\cup Z_{i}$ and let 
\begin{equation*}
A(Y_{1},Y_{2},\ldots ,Y_{k},Z_{1},Z_{2},\ldots ,Z_{k})=\left\{
(g_{1},g_{2},\ldots ,g_{k})\mid Y_{i}^{g_{i}}=Z_{i}\text{ for }1\leq i\leq
k\right\} \text{.}
\end{equation*}

Let us denote $S_{i}=G_{Y_{i}}$ the pointwise stabilizer of $Y_{i}$ in $G$,
let $C_{i}=S_{i}h_{i}$ and let $S=S_{1}\times S_{2}\times \cdots \times
S_{k} $. Obviously, $A$ coincides the right $S$-coset $C_{1}\times
C_{2}\times \cdots \times C_{k}$ so it inherits the natural topology and
measure from $S$.

Now let $x_{0}\in X\backslash Y$ be a fixed point and let $w\in F_{k}$ be a
reduced word of length $n>0$. Call a tuple $(g_{1},g_{2},\ldots ,g_{k})\in A$
\emph{distinctive}, if for the elements 
\begin{equation*}
x_{j}=x_{0}{}^{w_{j}(g_{1},g_{2},\ldots ,g_{n})}
\end{equation*}
we have $x_{i}\notin Y\cup Z$ ($1\leq i\leq n$) and $x_{i}\neq x_{j}$ ($%
0\leq i<j\leq n$). Let 
\begin{equation*}
U=U(Y_{1},Y_{2},\ldots ,Y_{k},Z_{1},Z_{2},\ldots ,Z_{k},x_{0},w)
\end{equation*}
be the set of distinctive elements of $A$.

Our first claim is that $U$ is open in $A$. Let $H=G_{x_{0}}$ be the
stabilizer of $x_{0}$ in $G$. For all $0\leq i<j\leq k$ let us define 
\begin{equation*}
V_{ij}=\left\{ (g_{1},g_{2},\ldots ,g_{k})\in A\mid x_{i}=x_{j}\right\} 
\text{. }
\end{equation*}
Clearly, 
\begin{equation*}
U=A\backslash \bigcup_{i<j\leq n}V_{ij}
\end{equation*}
so it is enough to show that the sets $V_{ij}$ are closed. Define the
function $F_{ij}:G^{k}\rightarrow G$ by 
\begin{equation*}
F_{ij}(g_{1},\ldots ,g_{k})=w_{i}(g_{1},g_{2},\ldots
,g_{n})w_{j}^{-1}(g_{1},g_{2},\ldots ,g_{n})\text{.}
\end{equation*}
Since $F_{ij}$ is continuous, the pre-image $W_{ij}=F_{ij}^{-1}(H)$ is
closed in $G^{n}$. On the other hand, $V_{ij}=A\cap W_{ij}$ so it is closed.

Our second claim is that $U$ is a.s. in $A$. We use induction on $n$, the
length of $w$. If $n=0$, then $U=A$ and the claim obviously holds.

If $n>0$, let $f$ be the first letter of $w$. By permuting the letters and
possibly switching $f_{1}$ and $f_{1}^{-1}$ (and also $Y_{1}$ and $Z_{1}$),
we can assume that $f=f_{1}$. Now let $R$ be the stabilizer of $x_{0}$ in $%
S_{1}$ and let $D=Rh_{1}\times C_{2}\times \cdots \times C_{k}$. Since $G$
acts strongly topologically on $X$, $R$ is closed of zero measure in $S_{1}$%
, thus $D$ is closed of zero measure in $A$. The set 
\begin{equation*}
L=\left\{ (g_{1},g_{2},\ldots ,g_{k})\in A\mid x_{0}^{g_{1}}\in Y\cup
\{x_{0}\}\right\}
\end{equation*}
is the finite union of $D$-cosets so it has zero measure in $A$. That is,
for almost all $(g_{1},g_{2},\ldots ,g_{k})\in A$ we have $%
x_{0}^{g_{1}}\notin Y\cup \{x_{0}\}$. For a fixed $x_{0}^{\prime
}=x_{0}^{g_{1}}$ consider the new system 
\begin{equation*}
Y_{1}^{\prime }=Y_{1}\cup \{x_{0}\}\text{, }Z_{1}^{\prime }=Z_{1}\cup
\{x_{0}^{\prime }\}
\end{equation*}
and $Y_{j}^{\prime }=Y_{j}$, $Z_{j}^{\prime }=Z_{j}$ for $2\leq j\leq k$.
Also let $w^{\prime }=f_{1}^{-1}w$ reduced (delete the first letter of $w$).
Now $w^{\prime }$ has length $n-1$ so by induction $U^{\prime
}=U(Y_{1}^{\prime },\ldots ,Z_{k}^{\prime },x_{0}^{\prime },w^{\prime })$ is
a.s. in $A^{\prime }=A(Y_{1}^{\prime },\ldots ,Z_{k}^{\prime })$. Since $%
x_{0}\in Y_{1}^{\prime }$, distinctive elements of $A^{\prime }$ will be
also distinctive in $A$. Using Lemma \ref{izelemma} we see that $U$ is a.s.
in $A$ what we claimed.

Now we prove the theorem. Let $k\geq 1,$ set $X_{1},\ldots ,X_{k+1}$ and $%
Y_{1}\ldots ,Y_{k+1}$ to be all empty sets and let $x_{0}\in X$. For a word $%
1\neq v\in F_{k}$ let us define $w\in F_{k+1}$ by $w=f_{k+1}vf_{k+1}^{-1}$.
Then using the above claim we get that the set of distinctive $k+1$-tuples
is a.s. in $G^{k+1}$. In particular 
\begin{equation*}
(x_{0}^{g_{k+1}})^{v(g_{1},\ldots ,g_{k})g^{-1}_{k+1}}=x_{0}^{w(g_{1},\ldots
,g_{k+1})}\neq x_{0}
\end{equation*}
so 
\begin{equation*}
(x_{0}^{g_{k+1}})^{v(g_{1},\ldots ,g_{k})} \neq x_{0}^{g_{k+1}}
\end{equation*}
for such $k+1$-tuples. Using Fubini's theorem this implies that for almost
all $h\in G$ the set of distinctive $k+1$-tuples with $g_{k+1}=h$ starting
from $x_{0}^{h}$ is a.s. in $G^{k}$. So assuming the transitivity of $G$ on $%
X$ for almost all $x\in X$ the value $v(g_{1},\ldots ,g_{k})$ does not fix $%
x $. Using that there are countably many words in $F_{k}$ and Fubini's
theorem again, we get that for almost all $k$-tuples $(g_{1},\ldots
,g_{k})\in G^{k}$ for almost all $x\in X$ the value $v(g_{1},\ldots ,g_{k})$
does not fix $x$ for all $1\neq v\in F_{k}$. That is, $(g_{1},\ldots
,g_{k})\in G^{k}$ generates a free subgroup acting almost freely on $X$ a.s. 
$\Box $

\bigskip

\textbf{Proof of Corollary \ref{regtree}. }Let $X$ be the boundary of the
tree, i.e., the set of infinite rays of $T$ and let $G=Aut(T)$. Then $G$ is
locally compact, the neighbourhood of the identity being pointwise
stabilizers of balls. We can assume that vertex stabilizers have Haar
measure $1$. Clearly, for every $x\in X$ the point stabilizer $G_{x}$ is
closed. Let $Y\subseteq X$ be an arbitrary (possibly empty) finite subset
and let $x\in X\backslash Y$.

If $\left| Y\right| \geq 3$, then $G_{Y}$ is compact as it the stabilizer of
the convex hull $C$ of $Y$ in $T$. Since $C$ has only finitely many limit
points in $X$ (the set $Y$), $x$ has a closest point $c$ in $C$. Then $G_{Y}$
acts transitively on the shadow of $c$ containing $x$, so the $G_{Y}$-orbit
containing $x$ is infinite (in fact it contains continuously many points).
As $G_{Y}$ is compact, this shows that $G_{Y\cup \{x\}}$ has infinite index
and thus zero measure in $G_{Y}$.

If $\left| Y\right| \leq 2$, then $G_{Y}$ acts transitively on $X\backslash
Y $. Fix a point $t\in T$, let $H$ be the stabilizer of $t$ in $G_{Y}$ and
let $H_{x}$ be the stabilizer of $x$ in $H$. Then $H$ is compact and has
countable index in $G_{Y}$ so it has positive measure in $G_{Y}$. As the
convex cone of $Y\cup \{t\}$ has no new boundary points, the $x$-orbit of $H$
is infinite and so $H_{x}$ has zero measure in $H$. Since $G_{Y\cup \{x\}}$
is a countable union of cosets of $H_{x}$, it has zero measure in $G_{Y}$.

The above argument shows that $G$ separates $X$ and the action is strongly
topological. Using Theorem \ref{topfree} we see that the corollary holds. $%
\Box $

\bigskip

Now we prove Theorem \ref{alter} using a general lemma on finite permutation
groups.

Let $G$ be a permutation group acting on a finite set $X$. We say that $G$
separates $X$ in order $(n,a)$ if for any subset $Y\subseteq X$ with $\left|
Y\right| =n$ the point stabilizer $G_{Y}$ acts on $X\backslash Y$ in a way
such that all orbits have size at least $a$. For a word $w\in F_{k}$ let $%
P_{G}(w)$ denote the probability that $w$ is not satisfied in $G$.

\begin{lemma}
\label{finperm}Let $G$ be a finite group separating $X$ of order $(n,a)$ and
let $w\in F_{k}$ have length $n$. Then $Q(w)\geq \left( 1-\frac{n}{a}\right)
^{n}$.
\end{lemma}

\textbf{Proof. }We will use the notations of the proof of Theorem \ref{nolaw}%
.

Let $x_{0}\in X$ a fixed point. We estimate $Q(w)$ from below with the
probability that $x_{0},x_{1},\ldots ,x_{n}$ are all disjoint. By induction
on $n$, with probability at least $\left( 1-\frac{n-1}{a}\right) ^{n-1}\geq
\left( 1-\frac{n}{a}\right) ^{n-1}$ the points $x_{0},x_{1},\ldots ,x_{n-1}$
are all disjoint. Fix such a $x_{0},x_{1},\ldots ,x_{n-1}$: by doing that,
we are restricting the probability space $G^{k}$ to the direct product of
particular cosets of stabilizers in $G$. Now the possible extensions for $%
x_{n}$ are uniformly chosen from a $G_{Y}$-orbit which does not meet $Y$.
Thus this orbit has size at least $a$ and the intersection of this orbit
with the set $x_{0},x_{1},\ldots ,x_{n-1}$ has size at most $n$. So with at
least probability $1-\frac{n}{a}$ the randomly chosen $x_{n}$ does not meet $%
x_{0},x_{1},\ldots ,x_{n-1}$, proving our lemma. $\Box $

\bigskip

\textbf{Proof of Theorem \ref{alter}. }Let $A=A_{k}$ be an alternating upper
composition factor of $G$. Using standard technique there exist $m$ and
normal subgroups $N\vartriangleleft M\vartriangleleft G$ such that $%
M/N\approxeq A^{m}$ and $\left| G:N\right| $ is finite. Factoring out with
the centralizer of $M/N$ in $G$ we get that $G$ quotients to a finite group $%
\widetilde{G}$ which embeds as a subgroup of the automorphism group of $%
A^{m} $ which contains $A^{m}$ and acts transitively on the factors.

Now the whole $Aut(A^{m})$ is isomorphic to the wreath product $S_{k}\wr
S_{m}$ which acts naturally on a set $X$ of $km$ points. Since the natural
action of $A$ on $k$ points separates it of order $(n,k-n)$ for any $n<k-1$,
the same holds for the action of $A^{m}$ on $X$. Now $\widetilde{G}$
contains $A^{m}$, so it also separates $X$ of order $(n,k-n)$ for all $n<k-1$%
.

Now let $w\in F_{k}$ be a word of length $n$. If a word is not satisfied in
the quotient group $\widetilde{G}$ then it is not satisfied for any preimage
in $G$ and the image of a uniform random element in $G$ is uniform random in 
$\widetilde{G}$. Using Lemma \ref{finperm}, we have 
\begin{equation*}
P_{G}(w)\geq P_{\widetilde{G}}(w)\geq \left( 1-\frac{n}{k-n}\right) ^{n}%
\text{.}
\end{equation*}
Since $k$ is unbounded, this implies $P_{G}(w)=1$.

Since there are countably many words in $F_{k}$, it follows that with
probability $1$ random $k$-tuples does not satisfy any word, that is, they
generate free subgroups. $\Box $

\bigskip

\noindent \textbf{Acknowledgement.} The author is grateful to the `Secret
Seminar' for the inspiring milieu and to Laci Pyber, who asked him whether
Theorem \ref{alter} holds.

\end{document}